\documentclass{elsart}

\usepackage{epsf,amssymb,latexsym,amsmath}
\usepackage{graphicx}
\usepackage{tikz}
\usepackage{mathabx}

\date{}
\def\NN{\hbox{\sf I\kern-.13em\hbox{N}}}
\def\RR{\hbox{\sf I\kern-.14em\hbox{R}}}
\def\ZZ{{\hbox{\sf Z\kern-.43emZ}}}

\def\Aut {\mathop{\rm Aut}\nolimits}

\def\mod {\mathop{\rm mod}\nolimits}
\def\s{\scriptstyle }

\newtheorem{theorem}{Theorem}[section]

\newtheorem{lemma}[theorem]{Lemma}

\newtheorem{corollary}[theorem]{Corollary}
\newtheorem{example}[theorem]{Example}

\begin{document}
 
 \begin{frontmatter}%********************************
 
\title{Optimal $L(2,1)$-labeling of certain strong graph bundles cycles over cycles }

\author{Irena Hrastnik Ladinek} 
\ead{irena.hrastnik@um.si}

\address{
FME,
University of Maribor,\\ Smetanova 17, Maribor 2000, Slovenia.\\} 

\begin{abstract} 
An $L(2,1)$-labeling of a graph $G=(V,E)$ is a function $f$ from the vertex set $V(G)$ to the set of nonnegative integers such that the labels on adjacent vertices differ by at least two, and the labels on vertices at  distance  two differ by at least one. The span of $f$ is the difference between the largest and the smallest numbers of $f(V)$. The  $\lambda$-number of $G$, denoted by $\lambda (G)$, is the minimum span over all $L(2,1)$-labelings of $G$. 
We prove that if  $X= C_m\boxtimes^{\sigma_\ell} C_{n}$ is a direct graph bundle with fiber $C_{n}$ and base $C_m$, $n$ is a multiple of 11 and $\ell$ has a form of $\ell =[11k+(-1)^a 4m]\mod n$ or of $\ell =[11k+(-1)^a 3m]\mod n$, where $a\in \{1,2\}$ and $k\in \ZZ$, then $\lambda (X)=10$.
\end{abstract}


 \end{frontmatter} %********************************

\section{Introduction}

The Frequency Assignment Problem (FAP) is a combinatorial optimization problem that arises in the field of telecommunications and radio frequency (RF) engineering. The goal of the FAP is to assign a set of communication frequencies to a set of transmitters while satisfying certain constraints and minimizing interference. This problem was first formulated as a graph coloring problem in 1980 by Hale \cite{Hale}. By Roberts \cite{Roberts}, in order to avoid interference, any two ''close'' transitters must receive different channels and any two ''very close'' transmitters must receive channels that are at least two channels apart. To translate the problem into the language of graph theory, the transmitters are represented by the vertices of a graph; two vertices are ''very close'' if they are adjacent and ''close'' if they are of distance two in the graph. Based on this problem,  Griggs and Yeh \cite{GrYe} introduced $L(2,1)$-labeling on a simple graph. 

Formally, an $L(2,1)$-labeling of a graph $G$ is an assignment $f$ of non-negative integers to vertices of $G$ such that 
$$|f(u)-f(v)|\geq \left\{\begin{array}{lr}2;& d(u,v)=1,\\
1;& d(u,v)=2.
\end{array}
\right.$$ 
We call these requirements the $L(2,1)$-{\it conditions} and $f(v)$ the {\it label} of $v$ under $f$. The $\lambda$-number of $G$, denoted by $\lambda(G)$, is the minimum value $\lambda$ such that $G$ admits an $L(2,1)$-labeling. 

Griggs and Yeh \cite{GrYe} put forward a conjecture that $\lambda(G) \leq \Delta^2$ for any graph $G$ with maximum degree $\Delta \geq 2$. They also proved that $\lambda(G)\leq\Delta^2+2\Delta$ for a general graph with maximum degree $\Delta$. Later, Chang and Kuo \cite{ChaKuo} improved the bound to $\Delta^2+\Delta$ while  Kr\'al' and \v Skrekovski \cite{Kral}  further reduced the bound to $\Delta^2+\Delta-1$. Moreover, in 2008, Gon\c calves \cite{Gonc} announced the bound $\Delta^2+\Delta-2$, the conjecture is still open. 

There are also a number of studies on the algorithms for the $L(2,1)$- labeling problem \cite {Bod, Fiala, Krat}. It is known to be ${\bf NP}$-hard for general graphs \cite{GeoMau}. Even for some relatively simple families of graphs such as planar graphs, bipartite graphs, chordal graphs \cite{Bod}, and graphs of treewidth two \cite{FialaGo}, the problem is also ${\bf NP}$-hard. Until now, only a few graph classes such as paths, cycles, and wheels are known to have polynomial  time algorithms for this problem \cite{Bod}.

Graph products are one of the natural constructions giving more complex graphs from simple ones.
Graph bundles \cite{Pisanski,PisanskiVra} generalize the notion of covering graphs and Cartesian products of graphs. The notion follows the definition of fiber bundles and vector bundles that became standard objects in topology \cite{Huse} as spaces which locally look like a product. Graph bundles corresponding to arbitrary graph products were introduced in \cite{PisanskiVra}.
While the labeling problem of graph products are well studied \cite{PKJha,Klav,KlaVes,KoVe,Shao,ShaoVe}, much less is known about the labeling problem of graph bundles \cite{KlavMo,Pisanski}.

In this paper, we focus on the $L(2,1)$-labeling problem of strong graph bundles cycles over cycles. We prove that the $\lambda$-number of strong graph bundle $C_m\boxtimes^{\sigma_\ell} C_n$ is 10, if certain conditions are imposed on the lengths of the cycles and on the  cyclic $\ell$-shift $\sigma_\ell$. After a section with basic terminology and some preliminary observations that are needed for the outline of our results, a detailed proof follows in Section 3, and the result is illustrated with an example.

\section{Terminology and Preliminaries} 

A finite, simple and undirected graph $G=(V(G),E(G))$ is given by a set of vertices $V(G)$ and a set of edges $E(G)$. As usual, the edge $\{i,j\}\in E(G)$ is shortly denoted by $ij$. Although we are interested in undirected graphs here, the order of the vertices will sometimes be important, 
 for example if we want to assign  automorphisms to the edges of the base graph.
 In such a case, we assign two opposite arcs $\{(i,j), (j,i)\}$ to the edge $\{i,j\}$. 

Denote by $P_n$ the {\em path} on $n\geq 1$ distinct vertices $0,1,2,\ldots, n-1$ with the  edges $ij \in E(P_n)$ if $j=i+1, 0\leq i <n-1.$The {\em cycle} $C_n$ on $n$ vertices is defined by  $V(C_n) = \{0,1,\dots,n-1\}$ and $ij \in E(C_n)$ if $i=(j\pm 1) \mod n$.
For a graph $G=(V,E)$ the distance $d_G(u,v)$, or briefly $d(u,v)$, between the vertices $u$ and $v$ is defined as the number of edges on a shortest path  between $u$ and $v$ in $G$.
  
Two graphs $G$ and $H$ are called {\em isomorphic}, in symbols $G\simeq H$, if there exists a bijection $\varphi $ 
from $V(G)$ onto $V(H)$ that preserves adjacency and nonadjacency.
In other words, a bijection $\varphi : V(G) \to V(H)$  is an  {\em isomorphism} when:   $\varphi(i)\varphi(j)\in E(H)$ 
if and only if $ij \in E(G)$.  An isomorphism of a graph $G$ onto itself is called an {\em automorphism}.
The identity automorphism on $G$ will be denoted by $id_G$ or shortly $id$.

There are two types of automorphisms of a cycle: cyclic shift of the cycle by $\ell$ elements will be briefly called cyclic $\ell$-shift and reflections with one, two or no fixed points (depending on the parity of $n$). We will focus on the first type. A cyclic $\ell$-shift, denoted by $\sigma_\ell$, $0\leq \ell <n$, defined as $\sigma_\ell(i)=(i+\ell)\mod n$ for $i=0,1,\ldots,n-1$. As a special case, we have the identity ($\ell=0$).

The strong product $G\boxtimes H$ of the graphs $G=(V(G),E(G))$ and $H=(V(H),E(H))$ is defined as follows:
$V (G \boxtimes H ) = V(G) \times V(H)$\\
$\begin{array}{ll}
E(G\boxtimes H)=  \{(g_1,h_1)(g_2,h_2)~|~ &h_1=h_2 {\mathop{\rm \;and\;}}g_1g_2\in E(G) {\mathop{\rm \;or\;}} \\
&g_1=g_2 {\mathop{\rm \;and\;}}h_1h_2\in E(H)\} {\mathop{\rm \;or\;}} \\
&g_1g_2\in E(G) {\mathop{\rm \;and\;}} h_1h_2\in E(H) \}.
\end{array}$\\
The strong product is commutative and associative. For more facts about the strong product of graphs we refer to \cite{ImKl}. 

Let $G$ and $H$ be graphs and let $\Aut(H)$ be the set of automorphisms of $H$.
To any ordered pair of adjacent vertices $g_1,g_2 \in V(G)$  we will assign an automorphism of $H$. 
Formally, let $\sigma : V(G)\times V(G) \to \Aut(H)$.
For brevity, we will write $\sigma(g_1,g_2)  = \sigma_{g_1,g_2}$
and assume  that $\sigma_{g_2,g_1} = \sigma_{g_1,g_2}^{-1}$ for any $g_1,g_2\in V(G)$.

Now we construct the graph $X$ as follows.
The vertex set of $X$ is the Cartesian product of the vertex sets, $V(X) = V(G) \times V(H)$.
The edges of $X$ are given by the rule:
\begin{itemize}
\item for any $g_1g_2 \in E(G)$ and any $h_1h_2 \in E(H)$, the vertices $(g_1,h_1)$ and 
$(g_2,\sigma_{g_1,g_2}(h_2)) $ are adjacent  in $X$  
\item for any  $h_1h_2 \in E(H)$ and any $g \in V(G)$, the vertices $(g,h_1)$ and $(g,h_2) $ are adjacent in $X$ 
\item for any $g_1g_2 \in E(G)$ and any  $h \in V(H)$, the vertices $(g_1,h)$ and $(g_2,\sigma_{g_1,g_2}(h)) $ are adjacent in $X$. 
\end{itemize}
We call $X$ a {\em strong graph bundle}  with base   $G$ and fiber $H$ and write $X = G\boxtimes^\sigma H$. Clearly, if all $\sigma_{g_1g_2}$ are identity automorphisms, the strong graph bundle is isomorphic to 
the strong product  $G\boxtimes^\sigma H = G\boxtimes  H$.

A graph bundle over a cycle can always be constructed in such a way that all but at most one automorphism are identities.
Fixing $V(C_n)= \{0,1,2,\dots,n-1\}$, let us denote $\sigma_{n-1,0} = \alpha$, 
$\sigma_{i-1,i} = id$ for $i=1,2,\dots,n-1$ and write   
$ C_n\boxtimes^\alpha H \simeq C_n\boxtimes^\sigma H $. We will use this facts  frequently in this article.

A graph bundle $ C_n\boxtimes^\alpha H$ can also be represented as the graph obtained from the product $P_n\boxtimes H$ by adding  a copy of $K_2\boxtimes H$ between the vertex sets $\{n-1\}\times V(H)$ and $\{0\}\times V(H)$ such that if $V(K_2)=\{1,2\}$ and $(1,u)$ is adjacent to $(2,v)$ in $K_2\boxtimes H$, then $(n-1,u)$ and $(0,\alpha(v))$ are connected by an edge in $ C_n\boxtimes^\alpha H$. 

The following lemma contains a useful lower bound for $\lambda(G)$ \cite{GrYe}.

\begin{lemma}
\label{lemma 1}
Let $G$ be a graph with maximum degree $\Delta\geq 2$. If $G$ contains three vertices of degree $\Delta$ such that one of them is adjacent to the other two, then $\lambda(G)\geq \Delta+2$.
\end{lemma} 

In the sequel we need the following facts.

\begin{claim}
\label{claim 1}
If $a,b$ and $n$ are integers with $n\geq 1$, then $|(a\mod n)-(b \mod n)|=(|a-b|\mod n)$ or $n-(|a-b|\mod n)$.
\end{claim}

\begin{corollary}
\label{posledica 1}
If $a,b,n$ and $p$ are integers with $n,p\geq 1$, then $|(a\mod n)-(b \mod n)|\geq p\Leftrightarrow p\leq  (|a-b|\mod n) \leq n-p$. 
\end{corollary}

\begin{corollary}
\label{posledica 2}
If $a,b$ and $n$ are integers with $n\geq 1$, then $|an-b| \mod n=(b \mod n)$ or $n-(b \mod n)$.
\end{corollary}

%%%%%%%%%%%%%%%%%%%%%%%%%%%%%%%%%%%%%%%%%%%%%%%%%%%%%%%%%%%%%%%%%%%%%%%%%%%%%%%%%%%%%%%%%%%%

\section{Result }

\begin{theorem}
\label{izrek1}
Let  $X= C_m\boxtimes^{\sigma_\ell} C_{n}, m\geq 3$ be a direct graph bundle with fiber $C_{n}$ and base $C_m$.  If $n$ is a multiple of 11 and $\ell$ has a form of $\ell =[11k+(-1)^a 4m]\mod n$ or of $\ell =[11k+(-1)^a 3m]\mod n$, where $a\in \{1,2\}$ and $k\in \ZZ$, then $\lambda(X)=10$.
\end{theorem}

\noindent {\bf Proof:}

To prove this theorem, we present four $L(2,1)$-labelings of $X$ using labels $0,1,\ldots , 10$ according to the cyclic $\ell$-shift $\sigma_\ell$. Let $v=(i,j)\in V(X)$.
\begin{enumerate}  
\item Let $\ell =[11k+(-1)^a4m]\mod n$, where $a\in \{1,2\}$ and $k\in \ZZ$. Define   labeling $f_a$ of $v$ as $$f_a(v)=\left[2i+(4+a)j\right]\mod 11$$ 
\item Let $\ell =[11k+(-1)^a3m]\mod n$, where $a\in \{1,2\}$ and $k\in \ZZ$. Define   labeling $g_a$ of $v$ as $$g_a(v)=\left[(4+a)i+2j\right]\mod 11$$
\end{enumerate}

All assignments are clearly well-defined. 

We will prove the theorem in two steps.
First, we consider a  spanning subgraph $P_m\boxtimes C_n$ of $ X=C_m\boxtimes^{\sigma_\ell} C_n$ in which the edges corresponding to the only 
(possibly) nontrivial automorphism are missing: edges between the vertex sets $\{m-1\}\times V(C_n)$ and $\{0\}\times V(C_n)$. These are edges  $(m-1,u)(0,(u+j'+\ell)\mod n)$, where $j'\in\{-1,0,1\}.$ Recall that $\sigma_{\ell}(i)=(i+\ell)\mod n$ and $\sigma_{\ell}^{-1}(i)=(i-\ell)\mod n$. Therefore, on the other hand, these are edges  $(0,u)(m-1,(u+j'-\ell)\mod n)$, $j'\in\{-1,0,1\}$.
 We will deal with these edges in the second step.

\begin{enumerate}
\item Consider the product $P_m\boxtimes C_n$. It is clear that $V(P_m\boxtimes C_n)=V(X)$. Let $X$ have labeling $f_a$. Then $v$ is assigned the integer  $$f_a(v)=\left[2i+(4+a)j\right]\mod 11.$$ 
Let $w$ be a vertex adjacent to $v$ in $P_m\boxtimes C_n$. Then $w$ is of the form $(i+i', (j+j')\mod n)$ with one of the following properties:
\begin{enumerate}
\item  if $i=0$, then $i'=1,j'\in\{-1,0,1\}$ or  $i'=0, j'\in\{-1, 1\}$,
\item  if $i=m-1$, then $i'=-1,j'\in\{-1, 0,1\}$ or $i'=0,j'\in\{-1, 1\}$, 
\item if $i=1,2,\ldots m-2$, then $i',j'\in\{-1, 0,1\}$ and $i',j'$ are not all zero.
\end{enumerate}

It is clear that $$f_a(w)=\left[2i+(4+a)j+2i'+(4+a)j'\right]\mod 11.$$
We will show that $v$ and $w$ receive labels that differ by at least two. By Corollary \ref{posledica 1}, to show that $|f_a(v)-f_a(w)|\geq 2$, it is enough to show that $$2\leq \left[|2i'+(4+a)j'|\mod 11\right] \leq 11-2=9$$  and, by symmetry, if $X$ has labeling $g_a$  that $$2\leq \left[|(4+a)i'+2j'|\mod 11\right] \leq 9.$$ The reader can easily verify that $|2i'+(4+a)j'|\mod 11,|(4+a)i'+2j'|\mod 11\in \{2,3,4,\ldots ,8\}$.

Now let $z$ be a vertex at a distance of two from $v$ in $P_m\boxtimes C_n$. Then $z$ is of the form $(i+i'', (j+j'')\mod n)$,  where $j''\in \{-2, -1, 0,1, 2\}$, $i''$ and $j''$ are not both different from 2 or -2, and one of the following properties holds:

\begin{enumerate}
\item $i=0, i''\in \{0,1,2\}$,
\item $i=1, i''\in \{-1,0,1,2\}$,
\item $i=m-1, i''\in \{-2,-1,0\}$,
\item $i=m-2, i''\in \{-2,-1,0,1\}$,
\item $i\in \{2,3,\ldots ,m-3\}, i''\in \{-2, -1, 0,1, 2\}$.
 
\end{enumerate}

Let $X$ have labeling $f_a$. Note that $z$ receives the label $$f_a(z)=\left[2i+(4+a)j+2i''+(4+a)j''\right]\mod 11.$$ We will show that $v$ and $z$ have different labels.
To show that $|f_a(v)-f_a(z)|\geq 1$, it is enough to show, by Corollary \ref{posledica 1}, that $$1\leq \left[|2i''+(4+a)j''|\mod 11\right] \leq 10 $$ and, 
if $X$ has  labeling $g_a$, that $$1\leq \left[|(4+a)i''+2j''|\mod 11\right] \leq 10 .$$
Since  $|2i''+(4+a)j''|\mod 11,|(4+a)i''+2j''|\mod 11\in \{1,2, 3, \ldots, 10\}\setminus \{7\}$, the calculation is simple and we omit it, it follows. 

\item Now we consider edges in $X= C_m\boxtimes^{\sigma_\ell} C_n$ between the fiber over $m-1$ and the fiber over $0$. 

\begin{enumerate}
\item Let us first consider two adjacent vertices, one from the fiber over $m-1$ and the other from the fiber over $0$. 

\begin{enumerate}
\item Suppose that $X$ has labeling $f_a$. Since $\ell =[11k+(-1)^a4m]\mod n$ for some $k\in \ZZ$, there exists a $k'\in \ZZ$ such that $\ell =11k'+(-1)^a 4m.$  

Let $v=(m-1,j)$ and let $w$ be of the form $w=(0,(j+j'+\ell)\mod n)$, where $j'\in \{-1,0,1\}$.
Then  $$f_a(v)=\left[2(m-1)+(4+a)j\right]\mod 11$$ and $$f_a(w)=\left[(4+a)(j+j'+\ell)\right]\mod 11.$$ To show that $|f_a(v)-f_a(w)|\geq 2$, it is enough to show that 
\begin{equation}
\label{a}
2\leq [|2(m-1)-(4+a)(j'+\ell)|\mod 11]\leq 9,
\end{equation} by Corollary \ref{posledica 1}.

Note that $$|2(m-1)-(4+a)(j'+\ell)|\mod 11=$$ $$|2(m-1)-(4+a)(j'+11k'+(-1)^a4m)|\mod 11=$$ 
$$|2m(1-(-1)^a(8+2a))-(2+(4+a)j')-(4+a)\cdot 11k')|\mod 11. $$
%$$=|2m\cdot11(-1)^{a+1}-(2+(4+a)j')-(4+a)\cdot11k')|\mod 11 $$

Since $1-(-1)^a(8+2a))$ is equal to $11$ for $a=1$ and equal to $-11$ for $a=2$, this is  equal to $$|2m\cdot11(-1)^{a+1}-(2+(4+a)j')-(4+a)\cdot11k')|\mod 11. $$

Assertion (\ref{a}) is true since, according to Corollary \ref{posledica 2},  $(2+(4+a)j')\mod 11\in \{2,7,8\}$ and $11-[(2+(4+a)j')\mod 11]\in  \{3,4,9\}$.

%Now let $w$ be of the form $w=(0, (j+\ell)\mod n)$. In this case, the vertex $w$ is assigned the integer  $$f_a(w)=\left[(4+a)(j+\ell)\right]\mod 11.$$  
%
%We claim that 
%$$2\leq [|2(m-1)-(4+a)\ell|\mod 11]\leq 9,$$ by Corollary \ref{posledica 1}. Note that $$|2(m-1)-(4+a)\ell|\mod 11=|2(m-1)-(4+a)(11k'+(-1)^a 4m)|\mod 11=$$ $$=|2m(1-(-1)^a(8+2a))-2-(4+a)\cdot 11k'|\mod 11=$$ $$=|2m\cdot 11(-1)^{a+1}-2-(4+a)\cdot 11k'|\mod 11. $$  By Corollary \ref{posledica 2},  this is equal to $2$ or $11-2=9$ and the claim is true. 

\item Let us now assume that $X$ has labeling $g_a$. Then $\ell =11k'+(-1)^a3m$ for some $k'\in \ZZ$ and  adjacent vertices $v=(m-1,j)$ and $w=(0,(j+j'+\ell)\mod n),j'\in \{-1,0,1\}$ are assigned the integers 
$$g_a(v)=\left[(4+a)(m-1)+2j\right]\mod 11$$ and $$g_a(w)=\left[2(j+j'+\ell)\right]\mod 11.$$
We must show that $|g_a(v)-g_a(w)|\geq 2,$ or by Corollary \ref{posledica 1}, that
\begin{equation}
\label{b}
2\leq [|(4+a)(m-1)-2(j'+\ell)|\mod 11]\leq 9.
\end{equation}

Note that $$|(4+a)(m-1)-2(j'+\ell)|\mod 11=$$ $$|(4+a)(m-1)-2(j'+11k'+(-1)^a3m)|\mod 11=$$
$$|m(4+a-6(-1)^a)-(4+a+2j')-22k'|\mod 11.$$
Since $4+a-6(-1)^a$ is equal to $11$ for $a=1$ and equal to $0$ for $a=2$ this is also  equal to $(4+a+2j')\mod 11\in \{3,4,5,6,7,8\}$  or to $11-[(4+a+2j')\mod 11] \in\{3,4,5,6,7,8\}$, by Corollary \ref{posledica 2}. The desired follows.

%For adjacent vertices $v=(m-1,j)$ and $w=(0,(j+\ell)\mod n)$ we have to show that 
%$$2\leq [|(4+a)(m-1)-2\ell|\mod 11]\leq 9.$$
%
%We find that $$|(4+a)(m-1)-2\ell|\mod 11=|(4+a)(m-1)-2(11k'+(-1)^a3m)|\mod 11=$$ $$|m(4+a-6(-1)^a)-(4+a)-22k'|\mod 11.$$ In this case too, since $(4+a)\mod 11,11-[(4+a)\mod 11]\in \{5,6\},$ the desired result follows.

\end{enumerate}

\item Finally, consider  vertices at a distance of two in $X$ where the shortest paths between them (there are of course several paths of length 2) contain at least one edge between the fiber over $m-1$ and the fiber over $0$. If  a path contains  two such edges, then vertices come from the same fiber, and we have already dealt with this.
We are therefore interested in vertices that contain exactly one edge between the fiber over $m-1$ and the fiber over $0$. So,  if  $v$ is from the fiber over $m-1$, then $z$ can be from the fiber over $0$ or over $1$, and if  $v$ is from the fiber over $m-2$, then $z$ is from the fiber over $0$. We claim that such two vertices get different labels. Suppose that $v=(m-1,j)$.

\begin{enumerate}
\item First, let $z$ be  from the fiber over $0$. Then $z$ is of the form  $z=(0,(j+j'+\ell)\mod n)$, where $j'\in \{ -2,2\}$. For $j'\in \{-1,0,1\}$,  $v$ and $z$ are adjacent vertices for which we have  proved that the received labels differ by at least two. It  suffices to replace $j'\in \{-1,0,1\}$ by $j'\in \{ -2,2\}$ in this proof and we obtain:

\begin{itemize}
\item  $|f_a(v)-f_a(z)|\geq 1$ since $(2+(4+a)j')\mod 11\in \{1,3\}$ and $11-[(2+(4+a)j')\mod 11]\in \{8,10\}$ and 
\item $|g_a(v)-g_a(z)|\geq 1$ since $(4+a+2j')\mod 11\in \{1,2,9,10\}$ and $11-[(2+(4+a)j')\mod 11]\in \{1,2,9,10\}.$ 
\end{itemize}

\item Now let $z$ be from the fiber over $1$. Then $z$ is of the form $z=(1,(j+j'+\ell)\mod n)$, where $j'\in\{-2,-1, 0,1, 2\}$.  Let $X$ have labeling $f_a$. Since  $v$ receives the label $$f_a(v)=\left[2(m-1)+(4+a)j\right]\mod 11$$ and $z$ receives the label $$f_a(z)=\left[2+(4+a)(j+j'+\ell)\right]\mod 11$$ it is enough to show that $$1\leq [|2(m-1)-2-(4+a)(j'+\ell)|\mod 11]\leq 10,$$ by Corollary \ref{posledica 1}.
In the proof of (\ref{a}), we have shown that $2(m-1)-(4+a)(j'+\ell)=11k''-(2+(4+a)j')$ for some $k''\in \ZZ$. Consequently $$|2(m-1)-2-(4+a)(j'+\ell)|\mod 11=$$ $$|11k''-(2+(4+a)j')-2|\mod 11=|11k''-(4+(4+a)j')|\mod 11.$$ Since $(4+(4+a)j')\mod 11\in \{3,4,5,9,10\}$ and $11-[(4+(4+a)j')\mod 11]\in \{1,2,6,7,8\},$  the claim is true. 

If $X$ has labeling $g_a$, then $$g_a(v)=\left[(4+a)(m-1)+2j\right]\mod 11$$ and $$g_a(z)=\left[(4+a)+2(j+j'+\ell)\right]\mod 11. $$ In this case we have to show that $$1\leq [|(4+a)(m-2)-2(j'+\ell)|\mod 11]\leq 10.$$ In the proof of (\ref{b}), we have shown that $(4+a)(m-1)-2(j'+\ell)=11k''-(4+a+2j')$ for some $k''\in \ZZ$. From this we obtain $$|(4+a)(m-2)-2(j'+\ell)|\mod 11=$$ $$|11k''-(4+a +2j')-(4+a)|\mod 11=|11k''-2(4+a +j')|\mod 11.$$ 

Since $2(4+a+j')\mod 11\in \{1,3,5,6,8,10\}$ and $11-[2(4+a+j'))\mod 11]\in \{1,3,5,6,8,10\},$  the desired follows.
\end{enumerate}

It is easy to see that the case $v=(m-2,j)$ and  $z=(0,(j+j'+\ell)\mod n)$, where $j'\in\{-2,-1, 0,1, 2\}$ is analogous to the above case (the first coordinate of the vertices $v=(m-1,j)$ and  $z=(1,(j+j'+\ell)\mod n)$ is reduced by 1). 

\end{enumerate}
\end{enumerate}

We have  shown that if $n$ is a multiple of $11$ and $\ell =[11k+(-1)^a4m]\mod n$ or $\ell =[11k+(-1)^a3m]\mod n$ for some $a\in \{1,2\}$ and  $k\in \ZZ$, then  $\lambda(C_m\boxtimes^{\sigma_\ell} C_n)\leq 10$. Since $C_m\boxtimes^{\sigma_\ell} C_n$ is also a regular graph of degree 8, an application of Lemma \ref{lemma 1} to the above statement is that $\lambda(C_m\boxtimes^{\sigma_\ell} C_n)= 10$.

\rule{2mm}{2mm}

\begin{example}
The $L(2,1)$-labeling of  $C_{13}\boxtimes^{\sigma_\ell} C_{11}$ with the labels $0,1,2\ldots ,10$ for $\ell=3,5,6,8$, based on 
the above scheme, is shown in Figure \ref{slika1} using the labeling of  $P_{13}\boxtimes P_{11}$.
\end{example}

\begin{figure}[h]
\begin{center}
\begin{tikzpicture} [scale=0.5]

	{    		
\foreach \x in {1,2,3,4,5,6,7,8,9,10,11,12,13}
 \foreach \y in {1,2,3,4,5,6,7,8,9,10,11}
		{%tocke
       \filldraw (\x,\y) circle (2pt);}		

%direktni
\draw  (1,1)--(11,11); \draw  (2,1)--(12,11); \draw  (3,1)--(13,11); \draw  (4,1)--(13,10); \draw  (5,1)--(13,9); \draw  (6,1)--(13,8); \draw  (7,1)--(13,7); \draw  (8,1)--(13,6); \draw  (9,1)--(13,5); \draw  (10,1)--(13,4); \draw  (11,1)--(13,3); \draw  (12,1)--(13,2); \draw  (13,1)--(13,1);

\draw  (1,2)--(10,11); \draw  (1,3)--(9,11); \draw  (1,4)--(8,11); \draw  (1,5)--(7,11); \draw  (1,6)--(6,11); \draw  (1,7)--(5,11); \draw  (1,8)--(4,11); \draw  (1,9)--(3,11); \draw  (1,10)--(2,11);

\draw  (1,2)--(2,1); \draw  (1,3)--(3,1); \draw  (1,4)--(4,1); \draw  (1,5)--(5,1); \draw  (1,6)--(6,1); \draw  (1,7)--(7,1);
\draw  (1,8)--(8,1); \draw  (1,9)--(9,1); \draw  (1,10)--(10,1); \draw  (1,11)--(11,1);

\draw  (2,11)--(12,1); \draw  (3,11)--(13,1); \draw  (4,11)--(13,2); \draw  (5,11)--(13,3); \draw  (6,11)--(13,4); \draw  (7,11)--(13,5); \draw  (8,11)--(13,6);
\draw  (9,11)--(13,7); \draw  (10,11)--(13,8); \draw  (11,11)--(13,9); \draw  (12,11)--(13,10); 

%kartezicni
 \foreach \x in {1,2,3,4,5,6,7,8,9,10,11,12}
 \foreach \y in {1,2,3,4,5,6,7,8,9,10,11}
		{%vodoravni
       \draw  (\x+0.6,\y)--(\x+1,\y);}

 \foreach \x in {1,2,3,4,5,6,7,8,9,10,11,12,13}
%navpicni
{\draw  (\x,1)--(\x,11);}

\path node at (1.4,1) {$\s 0$}; \path node at (2.4,1) {$\s 2$}; \path node at (3.4,1) {$\s 4$}; \path node at (4.4,1) {$\s 6$}; \path node at (5.4,1) {$\s 8$};
\path node at (6.4,1) {$\s 10$}; \path node at (7.4,1) {$\s 1$}; \path node at (8.4,1) {$\s 3$}; \path node at (9.4,1) {$\s 5$}; \path node at (10.4,1) {$\s 7$}; 
\path node at (11.4,1) {$\s 9$}; \path node at (12.4,1) {$\s 0$}; \path node at (13.4,1) {$\s 2$}; 

\path node at (1.4,2) {$\s 5$}; \path node at (2.4,2) {$\s 7$}; \path node at (3.4,2) {$\s 9$}; \path node at (4.4,2) {$\s 0$}; \path node at (5.4,2) {$\s2$};
\path node at (6.4,2) {$\s4$}; \path node at (7.4,2) {$\s6$}; \path node at (8.4,2) {$\s8$}; \path node at (9.4,2) {$\s10$};\path node at (10.4,2) {$\s1$}; 
\path node at (11.4,2) {$\s3$}; \path node at (12.4,2) {$\s5$}; \path node at (13.4,2) {$\s7$}; 

\path node at (1.4,3) {$\s10$}; \path node at (2.4,3) {$\s1$}; \path node at (3.4,3) {$\s3$}; \path node at (4.4,3) {$\s5$}; \path node at (5.4,3) {$\s7$};
\path node at (6.4,3) {$\s9$}; \path node at (7.4,3) {$\s0$}; \path node at (8.4,3) {$\s2$}; \path node at (9.4,3) {$\s4$};\path node at (10.4,3) {$\s6$}; 
\path node at (11.4,3) {$\s8$}; \path node at (12.4,3) {$\s10$}; \path node at (13.4,3) {$\s1$}; 

\path node at (1.4,4) {$\s4$}; \path node at (2.4,4) {$\s6$}; \path node at (3.4,4) {$\s8$}; \path node at (4.4,4) {$\s10$}; \path node at (5.4,4) {$\s1$};
\path node at (6.4,4) {$\s3$}; \path node at (7.4,4) {$\s5$}; \path node at (8.4,4) {$\s7$}; \path node at (9.4,4) {$\s9$};\path node at (10.4,4) {$\s0$}; 
\path node at (11.4,4) {$\s2$}; \path node at (12.4,4) {$\s4$}; \path node at (13.4,4) {$\s6$}; 

\path node at (1.4,5) {$\s9$}; \path node at (2.4,5) {$\s0$}; \path node at (3.4,5) {$\s2$}; \path node at (4.4,5) {$\s4$}; \path node at (5.4,5) {$\s6$};
\path node at (6.4,5) {$\s8$}; \path node at (7.4,5) {$\s10$}; \path node at (8.4,5) {$\s1$}; \path node at (9.4,5) {$\s3$};\path node at (10.4,5) {$\s5$}; 
\path node at (11.4,5) {$\s7$}; \path node at (12.4,5) {$\s9$}; \path node at (13.4,5) {$\s0$}; 
 
\path node at (1.4,6) {$\s3$}; \path node at (2.4,6) {$\s5$}; \path node at (3.4,6) {$\s7$}; \path node at (4.4,6) {$\s9$}; \path node at (5.4,6) {$\s0$};
\path node at (6.4,6) {$\s2$}; \path node at (7.4,6) {$\s4$}; \path node at (8.4,6) {$\s6$}; \path node at (9.4,6) {$\s8$};\path node at (10.4,6) {$\s10$}; 
\path node at (11.4,6) {$\s1$}; \path node at (12.4,6) {$\s3$}; \path node at (13.4,6) {$\s5$}; 

\path node at (1.4,7) {$\s8$}; \path node at (2.4,7) {$\s10$}; \path node at (3.4,7) {$\s1$}; \path node at (4.4,7) {$\s3$}; \path node at (5.4,7) {$\s5$};
\path node at (6.4,7) {$\s7$}; \path node at (7.4,7) {$\s9$}; \path node at (8.4,7) {$\s0$}; \path node at (9.4,7) {$\s2$};\path node at (10.4,7) {$\s4$}; 
\path node at (11.4,7) {$\s6$}; \path node at (12.4,7) {$\s8$}; \path node at (13.4,7) {$\s10$}; 

\path node at (1.4,8) {$\s2$}; \path node at (2.4,8) {$\s4$}; \path node at (3.4,8) {$\s6$}; \path node at (4.4,8) {$\s8$}; \path node at (5.4,8) {$\s10$};
\path node at (6.4,8) {$\s1$}; \path node at (7.4,8) {$\s3$}; \path node at (8.4,8) {$\s5$}; \path node at (9.4,8) {$\s7$};\path node at (10.4,8) {$\s9$}; 
\path node at (11.4,8) {$\s0$}; \path node at (12.4,8) {$\s2$}; \path node at (13.4,8) {$\s4$}; 

\path node at (1.4,9) {$\s7$}; \path node at (2.4,9) {$\s9$}; \path node at (3.4,9) {$\s0$}; \path node at (4.4,9) {$\s2$}; \path node at (5.4,9) {$\s4$};
\path node at (6.4,9) {$\s6$}; \path node at (7.4,9) {$\s8$}; \path node at (8.4,9) {$\s10$}; \path node at (9.4,9) {$\s1$};\path node at (10.4,9) {$\s3$}; 
\path node at (11.4,9) {$\s5$}; \path node at (12.4,9) {$\s7$}; \path node at (13.4,9) {$\s9$}; 

\path node at (1.4,10) {$\s1$}; \path node at (2.4,10) {$\s3$}; \path node at (3.4,10) {$\s5$}; \path node at (4.4,10) {$\s7$}; \path node at (5.4,10) {$\s9$};
\path node at (6.4,10) {$\s0$}; \path node at (7.4,10) {$\s2$}; \path node at (8.4,10) {$\s4$}; \path node at (9.4,10) {$\s6$};\path node at (10.4,10) {$\s8$}; 
\path node at (11.4,10) {$\s10$}; \path node at (12.4,10) {$\s1$}; \path node at (13.4,10) {$\s3$}; 

\path node at (1.4,11) {$\s6$}; \path node at (2.4,11) {$\s8$}; \path node at (3.4,11) {$\s10$}; \path node at (4.4,11) {$\s1$}; \path node at (5.4,11) {$\s3$};
\path node at (6.4,11) {$\s5$}; \path node at (7.4,11) {$\s7$}; \path node at (8.4,11) {$\s9$}; \path node at (9.4,11) {$\s0$};\path node at (10.4,11) {$\s2$}; 
\path node at (11.4,11) {$\s4$}; \path node at (12.4,11) {$\s6$}; \path node at (13.4,11) {$\s8$}; 

\path node at (1,0) {$0$}; \path node at (2,0) {$1$}; \path node at (3,0) {$2$}; \path node at (4,0) {$3$}; \path node at (5,0) {$4$};
\path node at (6,0) {$5$}; \path node at (7,0) {$6$}; \path node at (8,0) {$7$};  \path node at (9,0) {$8$};
\path node at (10,0) {$9$}; \path node at (11,0) {$10$}; \path node at (12,0) {$11$};  \path node at (13,0) {$12$};

\path node at (0,1) {$0$}; \path node at (0,2) {$1$}; \path node at (0,3) {$2$}; \path node at (0,4) {$3$}; \path node at (0,5) {$4$};
\path node at (0,6) {$5$}; \path node at (0,7) {$6$}; \path node at (0,8) {$7$}; \path node at (0,9) {$8$}; \path node at (0,10) {$9$};
\path node at (0,11) {$10$}; 

\path node at (7,-1) {$a)$};

%%%%%%%%%%%%%%%%
%b)

\foreach \x in {15,16,17,18,19,20,21,22,23,24,25,26,27}
 \foreach \y in {1,2,3,4,5,6,7,8,9,10,11}
		{%tocke
       \filldraw (\x,\y) circle (2pt);}

\path node at (15.4,1) {$\s 0$}; \path node at (16.4,1) {$\s 2$}; \path node at (17.4,1) {$\s 4$}; \path node at (18.4,1) {$\s 6$}; \path node at (19.4,1) {$\s 8$};
\path node at (20.4,1) {$\s 10$}; \path node at (21.4,1) {$\s 1$}; \path node at (22.4,1) {$\s 3$}; \path node at (23.4,1) {$\s 5$}; \path node at (24.4,1) {$\s 7$}; 
\path node at (25.4,1) {$\s 9$}; \path node at (26.4,1) {$\s 0$}; \path node at (27.4,1) {$\s 2$}; 

\path node at (15.4,2) {$\s 6$}; \path node at (16.4,2) {$\s 8$}; \path node at (17.4,2) {$\s 10$}; \path node at (18.4,2) {$\s 1$}; \path node at (19.4,2) {$\s3$};
\path node at (20.4,2) {$\s5$}; \path node at (21.4,2) {$\s7$}; \path node at (22.4,2) {$\s9$}; \path node at (23.4,2) {$\s0$};\path node at (24.4,2) {$\s2$}; 
\path node at (25.4,2) {$\s4$}; \path node at (26.4,2) {$\s6$}; \path node at (27.4,2) {$\s8$}; 

\path node at (15.4,3) {$\s1$}; \path node at (16.4,3) {$\s3$}; \path node at (17.4,3) {$\s5$}; \path node at (18.4,3) {$\s7$}; \path node at (19.4,3) {$\s9$};
\path node at (20.4,3) {$\s0$}; \path node at (21.4,3) {$\s2$}; \path node at (22.4,3) {$\s4$}; \path node at (23.4,3) {$\s6$};\path node at (24.4,3) {$\s8$}; 
\path node at (25.4,3) {$\s10$}; \path node at (26.4,3) {$\s1$}; \path node at (27.4,3) {$\s3$}; 

\path node at (15.4,4) {$\s7$}; \path node at (16.4,4) {$\s9$}; \path node at (17.4,4) {$\s0$}; \path node at (18.4,4) {$\s2$}; \path node at (19.4,4) {$\s4$};
\path node at (20.4,4) {$\s6$}; \path node at (21.4,4) {$\s8$}; \path node at (22.4,4) {$\s10$}; \path node at (23.4,4) {$\s1$};\path node at (24.4,4) {$\s3$}; 
\path node at (25.4,4) {$\s5$}; \path node at (26.4,4) {$\s7$}; \path node at (27.4,4) {$\s9$}; 

\path node at (15.4,5) {$\s2$}; \path node at (16.4,5) {$\s4$}; \path node at (17.4,5) {$\s6$}; \path node at (18.4,5) {$\s8$}; \path node at (19.4,5) {$\s10$};
\path node at (20.4,5) {$\s1$}; \path node at (21.4,5) {$\s3$}; \path node at (22.4,5) {$\s5$}; \path node at (23.4,5) {$\s7$};\path node at (24.4,5) {$\s9$}; 
\path node at (25.4,5) {$\s0$}; \path node at (26.4,5) {$\s2$}; \path node at (27.4,5) {$\s4$}; 
 
\path node at (15.4,6) {$\s8$}; \path node at (16.4,6) {$\s10$}; \path node at (17.4,6) {$\s1$}; \path node at (18.4,6) {$\s3$}; \path node at (19.4,6) {$\s5$};
\path node at (20.4,6) {$\s7$}; \path node at (21.4,6) {$\s9$}; \path node at (22.4,6) {$\s0$}; \path node at (23.4,6) {$\s2$};\path node at (24.4,6) {$\s4$}; 
\path node at (25.4,6) {$\s6$}; \path node at (26.4,6) {$\s8$}; \path node at (27.4,6) {$\s10$}; 

\path node at (15.4,7) {$\s3$}; \path node at (16.4,7) {$\s5$}; \path node at (17.4,7) {$\s7$}; \path node at (18.4,7) {$\s9$}; \path node at (19.4,7) {$\s0$};
\path node at (20.4,7) {$\s2$}; \path node at (21.4,7) {$\s4$}; \path node at (22.4,7) {$\s6$}; \path node at (23.4,7) {$\s8$};\path node at (24.4,7) {$\s10$}; 
\path node at (25.4,7) {$\s1$}; \path node at (26.4,7) {$\s3$}; \path node at (27.4,7) {$\s5$}; 

\path node at (15.4,8) {$\s9$}; \path node at (16.4,8) {$\s0$}; \path node at (17.4,8) {$\s2$}; \path node at (18.4,8) {$\s4$}; \path node at (19.4,8) {$\s6$};
\path node at (20.4,8) {$\s8$}; \path node at (21.4,8) {$\s10$}; \path node at (22.4,8) {$\s1$}; \path node at (23.4,8) {$\s3$};\path node at (24.4,8) {$\s5$}; 
\path node at (25.4,8) {$\s7$}; \path node at (26.4,8) {$\s9$}; \path node at (27.4,8) {$\s0$}; 

\path node at (15.4,9) {$\s4$}; \path node at (16.4,9) {$\s6$}; \path node at (17.4,9) {$\s8$}; \path node at (18.4,9) {$\s10$}; \path node at (19.4,9) {$\s1$};
\path node at (20.4,9) {$\s3$}; \path node at (21.4,9) {$\s5$}; \path node at (22.4,9) {$\s7$}; \path node at (23.4,9) {$\s9$};\path node at (24.4,9) {$\s1$}; 
\path node at (25.4,9) {$\s2$}; \path node at (26.4,9) {$\s4$}; \path node at (27.4,9) {$\s6$}; 

\path node at (15.4,10) {$\s10$}; \path node at (16.4,10) {$\s1$}; \path node at (17.4,10) {$\s3$}; \path node at (18.4,10) {$\s5$}; \path node at (19.4,10) {$\s7$}; \path node at (20.4,10) {$\s9$}; \path node at (21.4,10) {$\s0$}; \path node at (22.4,10) {$\s2$}; \path node at (23.4,10) {$\s4$};\path node at (24.4,10) {$\s6$}; \path node at (25.4,10) {$\s8$}; \path node at (26.4,10) {$\s10$}; \path node at (27.4,10) {$\s1$}; 

\path node at (15.4,11) {$\s5$}; \path node at (16.4,11) {$\s7$}; \path node at (17.4,11) {$\s9$}; \path node at (18.4,11) {$\s0$}; \path node at (19.4,11) {$\s2$};
\path node at (20.4,11) {$\s4$}; \path node at (21.4,11) {$\s6$}; \path node at (22.4,11) {$\s8$}; \path node at (23.4,11) {$\s10$};\path node at (24.4,11) {$\s1$}; \path node at (25.4,11) {$\s3$}; \path node at (26.4,11) {$\s5$}; \path node at (27.4,11) {$\s7$}; 

\path node at (15,0) {$0$}; \path node at (16,0) {$1$}; \path node at (17,0) {$2$}; \path node at (18,0) {$3$}; \path node at (19,0) {$4$};
\path node at (20,0) {$5$}; \path node at (21,0) {$6$}; \path node at (22,0) {$7$};  \path node at (23,0) {$8$};
\path node at (24,0) {$9$}; \path node at (25,0) {$10$}; \path node at (26,0) {$11$};  \path node at (27,0) {$12$};

\path node at (21,-1) {$b)$};

%direktni
\draw  (15,1)--(25,11); \draw  (16,1)--(26,11); \draw  (17,1)--(27,11); \draw  (18,1)--(27,10); \draw  (19,1)--(27,9); \draw  (20,1)--(27,8); \draw  (21,1)--(27,7); \draw  (22,1)--(27,6); \draw  (23,1)--(27,5); \draw  (24,1)--(27,4); \draw  (25,1)--(27,3); \draw  (26,1)--(27,2); \draw  (27,1)--(27,1);

\draw  (15,2)--(24,11); \draw  (15,3)--(23,11); \draw  (15,4)--(22,11); \draw  (15,5)--(21,11); \draw  (15,6)--(20,11); \draw  (15,7)--(19,11); \draw  (15,8)--(18,11); \draw  (15,9)--(17,11); \draw  (15,10)--(16,11);

\draw  (15,2)--(16,1); \draw  (15,3)--(17,1); \draw  (15,4)--(18,1); \draw  (15,5)--(19,1); \draw  (15,6)--(20,1); \draw  (15,7)--(21,1);
\draw  (15,8)--(22,1); \draw  (15,9)--(23,1); \draw  (15,10)--(24,1); \draw  (15,11)--(25,1);

\draw  (16,11)--(26,1); \draw  (17,11)--(27,1); \draw  (18,11)--(27,2); \draw  (19,11)--(27,3); \draw  (20,11)--(27,4); \draw  (21,11)--(27,5); \draw  (22,11)--(27,6);
\draw  (23,11)--(27,7); \draw  (24,11)--(27,8); \draw  (25,11)--(27,9); \draw  (26,11)--(27,10); 

%kartezicni
\foreach \x in {15,16,17,18,19,20,21,22,23,24,25,26}
 \foreach \y in {1,2,3,4,5,6,7,8,9,10,11}
		{%vodoravni
       \draw  (\x+0.6,\y)--(\x+1,\y);}

 \foreach \x in {15,16,17,18,19,20,21,22,23,24,25,26,27}
%navpicni
{\draw  (\x,1)--(\x,11);}

%%%%%%%%%%%%%%%%%%%%%%%%%%%%%%%%%%%%%%%%%%%%%%%
%c)

\foreach \x in {1,2,3,4,5,6,7,8,9,10,11,12,13}
 \foreach \y in {-3,-4,-5,-6,-7,-8,-9,-10,-11,-12,-13}
		{%tocke
       \filldraw (\x,\y) circle (2pt);}	

%direktni
\draw  (1,-13)--(11,-3); \draw  (2,-13)--(12,-3); \draw  (3,-13)--(13,-3); \draw  (4,-13)--(13,-4); \draw  (5,-13)--(13,-5); \draw  (6,-13)--(13,-6); \draw  (7,-13)--(13,-7); \draw  (8,-13)--(13,-8); \draw  (9,-13)--(13,-9); \draw  (10,-13)--(13,-10); \draw  (11,-13)--(13,-11); \draw  (12,-13)--(13,-12); 

\draw  (1,-12)--(10,-3); \draw  (1,-11)--(9,-3); \draw  (1,-10)--(8,-3); \draw  (1,-9)--(7,-3); \draw  (1,-8)--(6,-3); \draw  (1,-7)--(5,-3); \draw  (1,-6)--(4,-3); \draw  (1,-5)--(3,-3); \draw  (1,-4)--(2,-3);

\draw  (1,-12)--(2,-13); \draw  (1,-11)--(3,-13); \draw  (1,-10)--(4,-13); \draw  (1,-9)--(5,-13); \draw  (1,-8)--(6,-13); \draw  (1,-7)--(7,-13);
\draw  (1,-6)--(8,-13); \draw  (1,-5)--(9,-13); \draw  (1,-4)--(10,-13); \draw  (1,-3)--(11,-13);

\draw  (2,-3)--(12,-13); \draw  (3,-3)--(13,-13); \draw  (4,-3)--(13,-12); \draw  (5,-3)--(13,-11); \draw  (6,-3)--(13,-10); \draw  (7,-3)--(13,-9); \draw  (8,-3)--(13,-8);
\draw  (9,-3)--(13,-7); \draw  (10,-3)--(13,-6); \draw  (11,-3)--(13,-5); \draw  (12,-3)--(13,-4); 

%kartezicni
 \foreach \x in {1,2,3,4,5,6,7,8,9,10,11,12}
 \foreach \y in {-3,-4,-5,-6,-7,-8,-9,-10,-11,-12,-13}
		{%vodoravni
       \draw  (\x+0.6,\y)--(\x+1,\y);}

 \foreach \x in {1,2,3,4,5,6,7,8,9,10,11,12,13}
%navpicni
{\draw  (\x,-3)--(\x,-13);}

\path node at (1.4,-13) {$\s 0$}; \path node at (2.4,-13) {$\s 5$}; \path node at (3.4,-13) {$\s 10$}; \path node at (4.4,-13) {$\s 4$}; \path node at (5.4,-13) {$\s 9$}; \path node at (6.4,-13) {$\s 3$}; \path node at (7.4,-13) {$\s 8$}; \path node at (8.4,-13) {$\s 2$}; \path node at (9.4,-13) {$\s 7$}; \path node at (10.4,-13) {$\s 1$}; \path node at (11.4,-13) {$\s 6$}; \path node at (12.4,-13) {$\s 0$}; \path node at (13.4,-13) {$\s 5$}; 

\path node at (1.4,-12) {$\s 2$}; \path node at (2.4,-12) {$\s 7$}; \path node at (3.4,-12) {$\s 1$}; \path node at (4.4,-12) {$\s 6$}; \path node at (5.4,-12) {$\s0$}; \path node at (6.4,-12) {$\s5$}; \path node at (7.4,-12) {$\s10$}; \path node at (8.4,-12) {$\s4$}; \path node at (9.4,-12) {$\s9$};\path node at (10.4,-12) {$\s3$}; \path node at (11.4,-12) {$\s8$}; \path node at (12.4,-12) {$\s2$}; \path node at (13.4,-12) {$\s7$}; 

\path node at (1.4,-11) {$\s4$}; \path node at (2.4,-11) {$\s9$}; \path node at (3.4,-11) {$\s3$}; \path node at (4.4,-11) {$\s8$}; \path node at (5.4,-11) {$\s2$}; \path node at (6.4,-11) {$\s7$}; \path node at (7.4,-11) {$\s1$}; \path node at (8.4,-11) {$\s6$}; \path node at (9.4,-11) {$\s0$};\path node at (10.4,-11) {$\s5$}; \path node at (11.4,-11) {$\s10$}; \path node at (12.4,-11) {$\s4$}; \path node at (13.4,-11) {$\s9$}; 

\path node at (1.4,-10) {$\s6$}; \path node at (2.4,-10) {$\s0$}; \path node at (3.4,-10) {$\s5$}; \path node at (4.4,-10) {$\s10$}; \path node at (5.4,-10) {$\s4$}; \path node at (6.4,-10) {$\s9$}; \path node at (7.4,-10) {$\s3$}; \path node at (8.4,-10) {$\s8$}; \path node at (9.4,-10) {$\s2$};\path node at (10.4,-10) {$\s7$}; \path node at (11.4,-10) {$\s1$}; \path node at (12.4,-10) {$\s6$}; \path node at (13.4,-10) {$\s0$}; 

\path node at (1.4,-9) {$\s8$}; \path node at (2.4,-9) {$\s2$}; \path node at (3.4,-9) {$\s7$}; \path node at (4.4,-9) {$\s1$}; \path node at (5.4,-9) {$\s6$};
\path node at (6.4,-9) {$\s0$}; \path node at (7.4,-9) {$\s5$}; \path node at (8.4,-9) {$\s10$}; \path node at (9.4,-9) {$\s4$};\path node at (10.4,-9) {$\s9$}; 
\path node at (11.4,-9) {$\s3$}; \path node at (12.4,-9) {$\s8$}; \path node at (13.4,-9) {$\s2$}; 
 
\path node at (1.4,-8) {$\s10$}; \path node at (2.4,-8) {$\s4$}; \path node at (3.4,-8) {$\s9$}; \path node at (4.4,-8) {$\s3$}; \path node at (5.4,-8) {$\s8$};
\path node at (6.4,-8) {$\s2$}; \path node at (7.4,-8) {$\s7$}; \path node at (8.4,-8) {$\s1$}; \path node at (9.4,-8) {$\s6$};\path node at (10.4,-8) {$\s0$}; 
\path node at (11.4,-8) {$\s5$}; \path node at (12.4,-8) {$\s10$}; \path node at (13.4,-8) {$\s4$}; 

\path node at (1.4,-7) {$\s1$}; \path node at (2.4,-7) {$\s6$}; \path node at (3.4,-7) {$\s0$}; \path node at (4.4,-7) {$\s5$}; \path node at (5.4,-7) {$\s10$};
\path node at (6.4,-7) {$\s4$}; \path node at (7.4,-7) {$\s9$}; \path node at (8.4,-7) {$\s3$}; \path node at (9.4,-7) {$\s8$};\path node at (10.4,-7) {$\s2$}; 
\path node at (11.4,-7) {$\s7$}; \path node at (12.4,-7) {$\s1$}; \path node at (13.4,-7) {$\s6$}; 

\path node at (1.4,-6) {$\s3$}; \path node at (2.4,-6) {$\s8$}; \path node at (3.4,-6) {$\s2$}; \path node at (4.4,-6) {$\s7$}; \path node at (5.4,-6) {$\s1$};
\path node at (6.4,-6) {$\s6$}; \path node at (7.4,-6) {$\s0$}; \path node at (8.4,-6) {$\s5$}; \path node at (9.4,-6) {$\s10$};\path node at (10.4,-6) {$\s4$}; 
\path node at (11.4,-6) {$\s9$}; \path node at (12.4,-6) {$\s3$}; \path node at (13.4,-6) {$\s8$}; 

\path node at (1.4,-5) {$\s5$}; \path node at (2.4,-5) {$\s10$}; \path node at (3.4,-5) {$\s4$}; \path node at (4.4,-5) {$\s9$}; \path node at (5.4,-5) {$\s3$};
\path node at (6.4,-5) {$\s8$}; \path node at (7.4,-5) {$\s2$}; \path node at (8.4,-5) {$\s7$}; \path node at (9.4,-5) {$\s1$};\path node at (10.4,-5) {$\s6$}; 
\path node at (11.4,-5) {$\s0$}; \path node at (12.4,-5) {$\s5$}; \path node at (13.4,-5) {$\s10$}; 

\path node at (1.4,-4) {$\s7$}; \path node at (2.4,-4) {$\s1$}; \path node at (3.4,-4) {$\s6$}; \path node at (4.4,-4) {$\s0$}; \path node at (5.4,-4) {$\s3$};
\path node at (6.4,-4) {$\s10$}; \path node at (7.4,-4) {$\s4$}; \path node at (8.4,-4) {$\s9$}; \path node at (9.4,-4) {$\s3$};\path node at (10.4,-4) {$\s8$}; 
\path node at (11.4,-4) {$\s2$}; \path node at (12.4,-4) {$\s7$}; \path node at (13.4,-4) {$\s1$}; 

\path node at (1.4,-3) {$\s9$}; \path node at (2.4,-3) {$\s3$}; \path node at (3.4,-3) {$\s8$}; \path node at (4.4,-3) {$\s2$}; \path node at (5.4,-3) {$\s7$};
\path node at (6.4,-3) {$\s1$}; \path node at (7.4,-3) {$\s6$}; \path node at (8.4,-3) {$\s0$}; \path node at (9.4,-3) {$\s5$};\path node at (10.4,-3) {$\s10$}; 
\path node at (11.4,-3) {$\s4$}; \path node at (12.4,-3) {$\s9$}; \path node at (13.4,-3) {$\s3$}; 

\path node at (1,-14) {$0$}; \path node at (2,-14) {$1$}; \path node at (3,-14) {$2$}; \path node at (4,-14) {$3$}; \path node at (5,-14) {$4$};
\path node at (6,-14) {$5$}; \path node at (7,-14) {$6$}; \path node at (8,-14) {$7$};  \path node at (9,-14) {$8$};
\path node at (10,-14) {$9$}; \path node at (11,-14) {$10$}; \path node at (12,-14) {$11$};  \path node at (13,-14) {$12$};

\path node at (0,-13) {$0$}; \path node at (0,-12) {$1$}; \path node at (0,-11) {$2$}; \path node at (0,-10) {$3$}; \path node at (0,-9) {$4$};
\path node at (0,-8) {$5$}; \path node at (0,-7) {$6$}; \path node at (0,-6) {$7$}; \path node at (0,-5) {$8$}; \path node at (0,-4) {$9$};
\path node at (0,-3) {$10$}; 

\path node at (7,-15) {$c)$};

%%%%%%%%%%%%%%%%%%%%%%%%%%%%
%d)

\foreach \x in {15,16,17,18,19,20,21,22,23,24,25,26,27}
 \foreach \y in {-3,-4,-5,-6,-7,-8,-9,-10,-11,-12,-13}
		{%tocke
       \filldraw (\x,\y) circle (2pt);}	

\path node at (15.4,-13) {$\s 0$}; \path node at (16.4,-13) {$\s 6$}; \path node at (17.4,-13) {$\s 1$}; \path node at (18.4,-13) {$\s 7$}; \path node at (19.4,-13) {$\s 2$}; \path node at (20.4,-13) {$\s 8$}; \path node at (21.4,-13) {$\s 3$}; \path node at (22.4,-13) {$\s 9$}; \path node at (23.4,-13) {$\s 4$}; \path node at (24.4,-13) {$\s 10$}; \path node at (25.4,-13) {$\s 5$}; \path node at (26.4,-13) {$\s 0$}; \path node at (27.4,-13) {$\s 6$}; 

\path node at (15.4,-12) {$\s 2$}; \path node at (16.4,-12) {$\s 8$}; \path node at (17.4,-12) {$\s 3$}; \path node at (18.4,-12) {$\s 9$}; \path node at (19.4,-12) {$\s4$}; \path node at (20.4,-12) {$\s10$}; \path node at (21.4,-12) {$\s5$}; \path node at (22.4,-12) {$\s0$}; \path node at (23.4,-12) {$\s6$};\path node at (24.4,-12) {$\s1$}; \path node at (25.4,-12) {$\s7$}; \path node at (26.4,-12) {$\s2$}; \path node at (27.4,-12) {$\s8$}; 

\path node at (15.4,-11) {$\s4$}; \path node at (16.4,-11) {$\s10$}; \path node at (17.4,-11) {$\s5$}; \path node at (18.4,-11) {$\s0$}; \path node at (19.4,-11) {$\s6$}; \path node at (20.4,-11) {$\s1$}; \path node at (21.4,-11) {$\s7$}; \path node at (22.4,-11) {$\s2$}; \path node at (23.4,-11) {$\s8$};\path node at (24.4,-11) {$\s3$}; \path node at (25.4,-11) {$\s9$}; \path node at (26.4,-11) {$\s4$}; \path node at (27.4,-11) {$\s10$}; 

\path node at (15.4,-10) {$\s6$}; \path node at (16.4,-10) {$\s1$}; \path node at (17.4,-10) {$\s7$}; \path node at (18.4,-10) {$\s2$}; \path node at (19.4,-10) {$\s8$}; \path node at (20.4,-10) {$\s3$}; \path node at (21.4,-10) {$\s9$}; \path node at (22.4,-10) {$\s4$}; \path node at (23.4,-10) {$\s10$};\path node at (24.4,-10) {$\s5$}; \path node at (25.4,-10) {$\s0$}; \path node at (26.4,-10) {$\s6$}; \path node at (27.4,-10) {$\s1$}; 

\path node at (15.4,-9) {$\s8$}; \path node at (16.4,-9) {$\s3$}; \path node at (17.4,-9) {$\s9$}; \path node at (18.4,-9) {$\s4$}; \path node at (19.4,-9) {$\s10$}; \path node at (20.4,-9) {$\s5$}; \path node at (21.4,-9) {$\s0$}; \path node at (22.4,-9) {$\s6$}; \path node at (23.4,-9) {$\s1$};\path node at (24.4,-9) {$\s7$}; \path node at (25.4,-9) {$\s2$}; \path node at (26.4,-9) {$\s8$}; \path node at (27.4,-9) {$\s3$}; 
 
\path node at (15.4,-8) {$\s10$}; \path node at (16.4,-8) {$\s5$}; \path node at (17.4,-8) {$\s0$}; \path node at (18.4,-8) {$\s6$}; \path node at (19.4,-8) {$\s1$}; \path node at (20.4,-8) {$\s7$}; \path node at (21.4,-8) {$\s2$}; \path node at (22.4,-8) {$\s8$}; \path node at (23.4,-8) {$\s3$};\path node at (24.4,-8) {$\s9$}; \path node at (25.4,-8) {$\s4$}; \path node at (26.4,-8) {$\s10$}; \path node at (27.4,-8) {$\s5$}; 

\path node at (15.4,-7) {$\s1$}; \path node at (16.4,-7) {$\s7$}; \path node at (17.4,-7) {$\s2$}; \path node at (18.4,-7) {$\s8$}; \path node at (19.4,-7) {$\s3$}; \path node at (20.4,-7) {$\s9$}; \path node at (21.4,-7) {$\s4$}; \path node at (22.4,-7) {$\s10$}; \path node at (23.4,-7) {$\s5$};\path node at (24.4,-7) {$\s0$}; \path node at (25.4,-7) {$\s6$}; \path node at (26.4,-7) {$\s1$}; \path node at (27.4,-7) {$\s7$}; 

\path node at (15.4,-6) {$\s3$}; \path node at (16.4,-6) {$\s9$}; \path node at (17.4,-6) {$\s4$}; \path node at (18.4,-6) {$\s10$}; \path node at (19.4,-6) {$\s5$}; \path node at (20.4,-6) {$\s0$}; \path node at (21.4,-6) {$\s6$}; \path node at (22.4,-6) {$\s1$}; \path node at (23.4,-6) {$\s7$};\path node at (24.4,-6) {$\s2$}; \path node at (25.4,-6) {$\s8$}; \path node at (26.4,-6) {$\s3$}; \path node at (27.4,-6) {$\s9$}; 

\path node at (15.4,-5) {$\s5$}; \path node at (16.4,-5) {$\s0$}; \path node at (17.4,-5) {$\s6$}; \path node at (18.4,-5) {$\s1$}; \path node at (19.4,-5) {$\s7$};
\path node at (20.4,-5) {$\s2$}; \path node at (21.4,-5) {$\s8$}; \path node at (22.4,-5) {$\s3$}; \path node at (23.4,-5) {$\s9$};\path node at (24.4,-5) {$\s4$}; \path node at (25.4,-5) {$\s10$}; \path node at (26.4,-5) {$\s5$}; \path node at (27.4,-5) {$\s0$}; 

\path node at (15.4,-4) {$\s7$}; \path node at (16.4,-4) {$\s2$}; \path node at (17.4,-4) {$\s8$}; \path node at (18.4,-4) {$\s3$}; \path node at (19.4,-4) {$\s9$};
\path node at (20.4,-4) {$\s4$}; \path node at (21.4,-4) {$\s10$}; \path node at (22.4,-4) {$\s5$}; \path node at (23.4,-4) {$\s0$};\path node at (24.4,-4) {$\s6$}; \path node at (25.4,-4) {$\s1$}; \path node at (26.4,-4) {$\s7$}; \path node at (27.4,-4) {$\s2$}; 

\path node at (15.4,-3) {$\s9$}; \path node at (16.4,-3) {$\s4$}; \path node at (17.4,-3) {$\s10$}; \path node at (18.4,-3) {$\s5$}; \path node at (19.4,-3) {$\s0$}; \path node at (20.4,-3) {$\s6$}; \path node at (21.4,-3) {$\s1$}; \path node at (22.4,-3) {$\s7$}; \path node at (23.4,-3) {$\s2$};\path node at (24.4,-3) {$\s8$}; \path node at (25.4,-3) {$\s3$}; \path node at (26.4,-3) {$\s9$}; \path node at (27.4,-3) {$\s4$}; 

\path node at (15,-14) {$0$}; \path node at (16,-14) {$1$}; \path node at (17,-14) {$2$}; \path node at (18,-14) {$3$}; \path node at (19,-14) {$4$};
\path node at (20,-14) {$5$}; \path node at (21,-14) {$6$}; \path node at (22,-14) {$7$};  \path node at (23,-14) {$8$};
\path node at (24,-14) {$9$}; \path node at (25,-14) {$10$}; \path node at (26,-14) {$11$};  \path node at (27,-14) {$12$};

%direktni
\draw  (15,-13)--(25,-3); \draw  (16,-13)--(26,-3); \draw  (17,-13)--(27,-3); \draw  (18,-13)--(27,-4); \draw  (19,-13)--(27,-5); \draw  (20,-13)--(27,-6); \draw  (21,-13)--(27,-7); \draw  (22,-13)--(27,-8); \draw  (23,-13)--(27,-9); \draw  (24,-13)--(27,-10); \draw  (25,-13)--(27,-11); \draw  (26,-13)--(27,-12); \draw  (27,-13)--(27,-13);

\draw  (15,-12)--(24,-3); \draw  (15,-11)--(23,-3); \draw  (15,-10)--(22,-3); \draw  (15,-9)--(21,-3); \draw  (15,-8)--(20,-3); \draw  (15,-7)--(19,-3); \draw  (15,-6)--(18,-3); \draw  (15,-5)--(17,-3); \draw  (15,-4)--(16,-3);

\draw  (15,-12)--(16,-13); \draw  (15,-11)--(17,-13); \draw  (15,-10)--(18,-13); \draw  (15,-9)--(19,-13); \draw  (15,-8)--(20,-13); \draw  (15,-7)--(21,-13);
\draw  (15,-6)--(22,-13); \draw  (15,-5)--(23,-13); \draw  (15,-4)--(24,-13); \draw  (15,-3)--(25,-13);

\draw  (16,-3)--(26,-13); \draw  (17,-3)--(27,-13); \draw  (18,-3)--(27,-12); \draw  (19,-3)--(27,-11); \draw  (20,-3)--(27,-10); \draw  (21,-3)--(27,-9); \draw  (22,-3)--(27,-8); \draw  (23,-3)--(27,-7); \draw  (24,-3)--(27,-6); \draw  (25,-3)--(27,-5); \draw  (26,-3)--(27,-4); 

%kartezicni
 \foreach \x in {15,16,17,18,19,20,21,22,23,24,25,26}
 \foreach \y in {-3,-4,-5,-6,-7,-8,-9,-10,-11,-12,-13}
		{%vodoravni
       \draw  (\x+0.6,\y)--(\x+1,\y);}

 \foreach \x in {15,16,17,18,19,20,21,22,23,24,25,26,27}
%navpicni
{\draw  (\x,-3)--(\x,-13);}

\path node at (21,-15) {$d)$};
}	
\end{tikzpicture}
\caption{Four $L(2,1)$-labelings of $P_{13}\times P_{11}$ towards that of  $C_{13}\times^{\sigma_\ell}C_{11}$corresponding to the cyclic $\ell$-shift $\sigma_\ell$: \\
a) $ \ell =3, f_1(i,j)=\left[2i+5j\right]\mod 11$,\\
b) $ \ell =8, f_2(i,j)=\left[2i+6j)\right]\mod 11 $, \\
c) $ \ell=5, g_1(i,j)=\left[5i+2j\right]\mod 11 $, \\
d) $ \ell=6, g_2(i,j)=\left[6i+2j\right]\mod 11$.}
\label{slika1}
\end{center}
\end{figure}

%%%%%%%%%%%%%%%%%%%%%%%%%%%%%%%%%%%%%%%%%%%%%%%%%%%%%%%%%
%%%%%%%%%%%%%%%%%%%%%%%%%%%%%%%%%%%%%%%%%%%%%%%%%%%%%%%%

\end{document}